\documentclass[runningheads,a4paper]{llncs}

\usepackage{subfigure}
\usepackage{graphicx}
\usepackage{amsmath}
\usepackage{amsfonts}
\usepackage{amssymb} 
\usepackage{color}
\usepackage{subfigure}
\usepackage{algorithm}
\usepackage{algpseudocode}

\DeclareMathOperator*{\argmin}{\arg \min}

\DeclareMathOperator{\jmc}{\mathcal{J  }}

\newcommand{\norm}[1]{\Vert #1 \Vert}%
\newcommand{\R}{\mathbb{R}}%

\newcommand{\jac}{\jmc\!}

\begin{document}

\mainmatter  

\title{Preconditioned ADMM with nonlinear\\ operator constraint}

\titlerunning{Preconditioned ADMM with nonlinear operator constraint}

%
%
\author{Martin Benning\inst{1} \and Florian Knoll\inst{2} \and Carola-Bibiane Sch\"{o}nlieb\inst{1}\and\\  Tuomo Valkonen\inst{1}}
\authorrunning{Benning et al.}

\institute{University of Cambridge, Department of Applied Mathematics and Theoretical Physics, Wilberforce Road, Cambridge CB3 0WA, United Kingdom\\ (\{mb941, cbs31, tjmv3\}@cam.ac.uk) \and New York University, Center for Advanced Imaging, Innovation and Research, New York 4th Floor 660 First Avenue New York, NY 10016, USA\\ 
(florian.knoll@nyumc.org)}

%
%

\toctitle{Preconditioned ADMM with nonlinear operator constraint}
\tocauthor{Martin Benning}
\maketitle

\begin{abstract}
We are presenting a modification of the well-known Alternating Direction Method of Multipliers (ADMM) algorithm with additional preconditioning that aims at solving convex optimisation problems with nonlinear operator constraints. Connections to the recently developed Nonlinear Primal-Dual Hybrid Gradient Method (NL-PDHGM) are presented, and the algorithm is demonstrated to handle the nonlinear inverse problem of parallel Magnetic Resonance Imaging (MRI).

 
\end{abstract}

\section{Introduction}
Non-smooth regularisation methods are popular tools in image processing. They allow to promote sparsity of inverse problem solutions with respect to specific representations; they allow to implicitly restrict the null-space of the forward operator while guaranteeing noise suppression at the same time. The most prominent representatives of this class are total variation regularisation \cite{rudin1992nonlinear} and $\ell^1$-norm regularisation as in the broader context of compressed sensing \cite{donoho2006compressed,candes2006compressive}.

In order to solve convex, non-smooth regularisation methods with linear operator constraints computationally, first-order operator splitting methods have gained increasing interest over the last decade, see \cite{gabay1983,goldstein2009split,beck2009fast,chambolle2011first} to name just a few. Despite some recent extensions to certain types of non-convex problems \cite{ochs2014ipiano,moellenhoff2015primal,Bonettini2015,7271087} there has to our knowledge only been made little progress for nonlinear operators constraints \cite{bachmayr2009iterative,valkonen2014primal}.

In this paper we are particularly interested in minimising non-smooth, convex functionals with nonlinear operator constraints. This model covers many interesting applications; one particular application that we are going to address is the joint reconstruction of the spin-proton-density and coil sensitivity maps in parallel MRI \cite{uecker2008image,knoll2012parallel}. 

The paper is structured as follows: we will introduce the generic problem formulation, then address its numerical minimisation via a generalised ADMM method with linearised operator constraints. Subsequently we will show connections to the recently proposed NL-PDHGM method (indicating a local convergence result of the proposed algorithm) and conclude with the joint spin-proton-density and coil sensitivity map estimation as a numerical example.

\section{Problem formulation}
We consider the following generic constrained minimisation problem:
\begin{align}
(\hat{u}, \hat{v}) &= \argmin_{u, v} \left\{ H(u) + J(v) \ \text{subject to} \ F(u, v) = c\right\} \, \text{.}\label{eq:nonlinprobform}
\end{align}
Here $H$ and $J$ denote proper, convex and lower semi-continuous functionals, $F$ is a nonlinear operator and $c$ a given function. Note that for nonlinear operators of the form $F(u, v) = G(u) - v$ and $c = 0$ problem \eqref{eq:nonlinprobform} can be written as 
\begin{align}
\hat{u} &= \argmin_{u} \left\{ H(u) + J(G(u)) \right\} \text{.}\label{eq:nonlinvarprob}
\end{align}

In the following we want to propose a strategy for solving \eqref{eq:nonlinprobform} that is based on simultaneous linearisation of the nonlinear operator constraint and the solution of an inexact ADMM problem. 

%

\section{Alternating direction method of multipliers}
We solve \eqref{eq:nonlinprobform} by alternating optimisation of the augmented Lagrange function
\begin{align}
\mathcal{L}_\delta(u, v; \mu) = H(u) + J(v) + \langle \mu, F(u, v) - c \rangle + \frac{\delta}{2} \| F(u, v) - c\|_2^2 \, \text{.}\label{eq:auglag}
\end{align}
Alternating minimisation of \eqref{eq:auglag} in 
$u$, $v$ and subsequent maximisation of $\mu$ via a step of gradient ascent yields this nonlinear version of ADMM \cite{gabay1983}:
\begin{align}
u^{k + 1} &\in \argmin_u \left\{ \frac{\delta}{2} \| F(u, v^k) - c\|_2^2 + \langle \mu^k, F(u, v^k) \rangle + H(u) \right\} \, \text{,} \label{eq:admmup1}\\
v^{k + 1} &\in \argmin_v \left\{ \frac{\delta}{2} \| F(u^{k + 1}, v) - c\|_2^2 + \langle \mu^k, F(u^{k + 1}, v) \rangle + J(v) \right\} \, \text{,} \label{eq:admmvp1} \\
\mu^{k + 1} &= \mu^k + \delta \left( F(u^{k + 1}, v^{k + 1}) - c \right) \, \text{.} \label{eq:admmupmu}
\end{align}
Not having to deal with nonlinear subproblems, we replace $F(u^{k + 1}, v^k)$ and $F(u^{k + 1}, v^{k + 1})$ by their Taylor linearisations around $u^k$ and $v^k$, which yields $F(u, v^k) \approx F(u^k, v^k) + \partial_u F(u^k, v^k)\left( u - u^k \right)$ and $F(u^{k + 1}, v) \approx F(u^{k + 1}, v^k) + \partial_v F(u^{k + 1}, v^k)\left( v - v^k \right)$, respectively. The updates \eqref{eq:admmup1} and \eqref{eq:admmvp1} modify to
\begin{align}
u^{k + 1} &\in \argmin_u \left\{ \frac{\delta}{2} \left\| A^k u - c_1^k\right\|_2^2 + \langle \mu^k, A^k u \rangle + H(u) \right\} \, \text{,} \label{eq:admmup2}\\
v^{k + 1} &\in \argmin_v \left\{ \frac{\delta}{2} \left\| B^k v - c_2^k \right\|_2^2 + \langle \mu^k, B^k v \rangle + J(v) \right\} \, \text{,} \label{eq:admmvp2}
\end{align}
with $A^k := \partial_u F(u^k, v^k)$, $B^k := \partial_v F(u^{k + 1}, v^k)$,  $c_1^k := c + A^k u^k - F(u^k, v^k)$ and $c_2^k := c + B^k v^k - F(u^{k + 1}, v^k)$. Note that the updates \eqref{eq:admmup2} and \eqref{eq:admmvp2} are still implicit, regardless of $H$ and $J$. In the following, we want to modify the updates such that they become simple proximity operations.

\section{Preconditioned ADMM}
Based on \cite{zhang2011unified}, we modify \eqref{eq:admmup2} and \eqref{eq:admmvp2} by adding the surrogate terms $\| u^{k + 1} - u^k \|_{Q^k_1}^2 / 2$ and $\| v^{k + 1} - v^k \|_{Q^k_2}^2 / 2$, with $\| w \|_Q := \sqrt{\langle Qw, w\rangle}$ (note that if $Q$ is chosen to be positive definite, $\| \cdot \|_Q$ becomes a norm). We then obtain
\begin{align*}
u^{k + 1} &\in \argmin_u \left\{ \frac{\delta}{2} \left\| A^k u - c_1^k\right\|_2^2 + \langle \mu^k, A^k u \rangle + H(u) + \frac{1}{2}\| u - u^k \|_{Q_1^k}^2 \right\} \, \text{,} \\
v^{k + 1} &\in \argmin_v \left\{ \frac{\delta}{2} \left\| B^k v - c_2^k \right\|_2^2 + \langle \mu^k, B^k v \rangle + J(v) + \frac{1}{2}\| v - v^k \|_{Q_2^k}^2 \right\} \, \text{.}
\end{align*}
If we choose $Q_1^k := \tau_1^k I - \delta A^k {}^* A^k$ with $\tau_1^k \delta < 1/\| A^k \|^2$ and $Q_2^k := \tau_2^k I - \delta B^k {}^* B^k$ with $\tau_2^k \delta < 1/\| B^k \|^2$ and if we define $\overline{\mu}^k := 2\mu^k - \mu^{k - 1}$ we obtain
\begin{align}
u^{k + 1} &= \left( I + \tau_1^k \partial H \right)^{-1} \left( u^k - \tau_1^k A^k{}^*\overline{\mu}^k  \right)\, \text{,} \label{eq:admmup3}\\
v^{k + 1} &= \left( I + \tau_2^k \partial J \right)^{-1} \left( v^k - \tau_2^k B^k{}^*\left( \mu^k + \delta \left(F(u^{k + 1}, v^k) - c \right) \right) \right) \, \text{,} \label{eq:admmvp3}
\end{align}
with $(I + \alpha \partial E)^{-1}(w)$ denoting the proximity or resolvent operator
\begin{align*}
(I + \alpha \partial E)^{-1}(w) := \argmin_u \left\{ \frac{1}{2} \|u - w\|_2^2 + \alpha E(u) \right\} \, \text{.}
\end{align*}
The entire proposed algorithm with updates \eqref{eq:admmup3}, \eqref{eq:admmvp3} and \eqref{eq:admmupmu} reads as
\begin{algorithm}[H]
	\caption{Preconditioned ADMM with nonlinear operator constraint}		
	\begin{algorithmic}
		\State \textbf{Parameters:} $H, \ J, \ F, \ c$
		\State \textbf{Initialization:} $u^0$, $v^0$, $\mu^0$, $\delta$
		\State $\overline{\mu}^0 = \mu^0$
			\While{convergence criterion is not met}
				\vspace{0.1cm}
				\State $A^k = \partial_u F(u^k, v^k)$
				\State Set $\tau_1^k$ such that $\tau_1^k \delta < 1/\| A^k \|^2$
				\State $u^{k + 1} = \left( I + \tau_1^k \partial H \right)^{-1} \left( u^k - \tau_1^k A^k{}^*\overline{\mu}^k  \right)$
				\State $B^k = \partial_v F(u^{k + 1}, v^k)$
				\State Set $\tau_2^k$ such that $\tau_2^k \delta < 1/\| B^k \|^2$
				\State $v^{k + 1} = \left( I + \tau_2^k \partial J \right)^{-1} \left( v^k - \tau_2^k B^k{}^*\left( \mu^k + \delta \left(F(u^{k + 1}, v^k) - c \right) \right) \right)$
				\State $\mu^{k + 1} = \mu^k + \delta \left( F(u^{k + 1}, v^{k + 1}) - c \right)$
				\State $\overline{\mu}^{k + 1} = 2\mu^{k + 1} - \mu^k$
				\vspace{0.1cm}
			\EndWhile
			\State \Return $u^{k}$, $v^k$, $\mu^k$, $\overline{\mu}^k$
	\end{algorithmic}
	\label{alg:nlpadmm}
\end{algorithm}

\section{Connection to NL-PDHGM}\label{sec:nlpdhgm}
In the following we want to show how the algorithm simplifies in case the nonlinear operator constraint is only nonlinear in one variable, which is sufficient for problems of the form \eqref{eq:nonlinvarprob}. Without loss of generality we consider constraints of the form $F(u, v) = G(u) - v$, where $G$ represents a nonlinear operator in $u$. Then we have $A^k = \jac G(u^k)$ (with $\jac G(u^k)$ denoting the Jacobi matrix of $G$ at $u^k$), $B^k = -I$ and if we further choose $\tau_2^k = 1/\delta$ for all $k$, update \eqref{eq:admmvp3} reads
\begin{align*}
v^{k + 1} = \left( I + \frac{1}{\delta} \partial J \right)^{-1} \left( G(u^{k + 1}) + \frac{1}{\delta}\mu^k \right) \, \text{.}
\end{align*}
Applying Moreau's identity \cite{rockafellar1970convex} $b = \left(I + \frac{1}{\delta} \partial J\right)^{-1}(b) + \frac{1}{\delta}(I + \delta \partial J^*)^{-1}(\delta b)$ yields 
\begin{align*}
\mu^{k + 1} = \left( I + \delta \partial J^{*} \right)^{-1}\left( \mu^k + \delta  G(u^{k + 1}) \right) \, \text{.}
\end{align*}
If we further change the order of the updates, starting with the update for $\mu$, the whole algorithm reads
\begin{align*}
\mu^{k + 1} &= \left( I + \delta \partial J^{*} \right)^{-1}\left( \mu^k + \delta G(u^k) \right) \, \text{,}\\
\overline{\mu}^{k + 1} &= 2\mu^{k + 1} - \mu^k \, \text{,}\\
u^{k + 1} &= \left( I + \tau_1^k \partial H \right)^{-1} \left( u^k - \tau_1^k \jac G(u^k)^* \overline{\mu}^{k + 1}  \right)\, \text{.}
\end{align*}
Note that this algorithm is almost the same as NL-PDHGM proposed in \cite{valkonen2014primal} for $\theta = 1$, except that the extrapolation step is carried out on the dual variable $\mu$ instead of the primal variable $u$. In the following we want to briefly sketch how to prove convergence for this algorithm in analogy to \cite{valkonen2014primal}. We define 
\begin{align*}
N(\mu^{k + 1}, u^{k + 1}) &:= \left( \begin{array}{c} 
\partial J^*(\mu^{k + 1}) - \nabla G(u^k) u^{k + 1} - c^k\\ 
\partial H(u^{k + 1}) + \jac G(u^k)^* \mu^{k + 1}
\end{array} \right) \, \text{,}\\
L^k &:= \left( \begin{array}{cc}
\frac{1}{\delta} I & \jac G(u^k)\\
\jac G(u^k)^* & \frac{1}{\tau_1^k} I
\end{array}\right) \, \text{,}
\end{align*}
with $c^k := G(u^k) - \jac G(u^k)u^k$.
Now the algorithm is: find $(\mu^{k+1}, u^{k+1})$ such that
\[
    N(\mu^{k + 1}, u^{k + 1}) + L^k(\mu^{k + 1}-\mu^k, u^{k + 1}-u^k) \ni 0.
\]  
If we exchange the order of $\mu$ and $u$ here, i.e., reorder the rows of $N$, and the rows and columns of $L^k$, we obtain \emph{almost} the ``linearised'' NL-PDHGM of \cite{valkonen2014primal}. The difference is that the sign of $\jmc G$ in $L^k$ is inverted. The only points in \cite{valkonen2014primal} where the exact structure of $L^k$ ($M_{x^k}$ therein) is used, are Lemma 3.1, Lemma 3.6 and Lemma 3.10. The first two go through exactly as before with the negated structure. Reproducing Lemma 3.10 demands bounding actual step lengths $\norm{u^k-u^{k+1}}$ and $\norm{\mu^k-\mu^{k+1}}$ from below, near a solution for arbitrary $\epsilon>0$. A proof would go beyond the page limit of this proceeding.  Let us just point out that this can be done, implying that the convergence results of \cite{valkonen2014primal} apply for this algorithm as well. This means that under somewhat technical regularity conditions, which for TV type problems amount to Huber regularisation, local convergence in a neighbourhood of the true solution can be guaranteed.

\section{Joint estimation of the spin-proton density and coil sensitivities in parallel MRI}
\begin{figure}[!ht]\label{fig:expsetup}
\begin{center}
\subfigure[Brain phantom]{\includegraphics[width=0.35\textwidth]{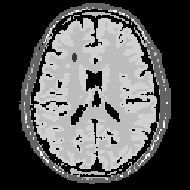}\label{subfig:gt1}}
\hspace{0.1cm}
\subfigure[25\% sub-sampling]{\includegraphics[width=0.35\textwidth]{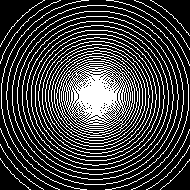}\label{subfig:sampscheme}\label{subfig:sampscheme}}\hspace{0.25cm}\includegraphics[width=0.0975\textwidth]{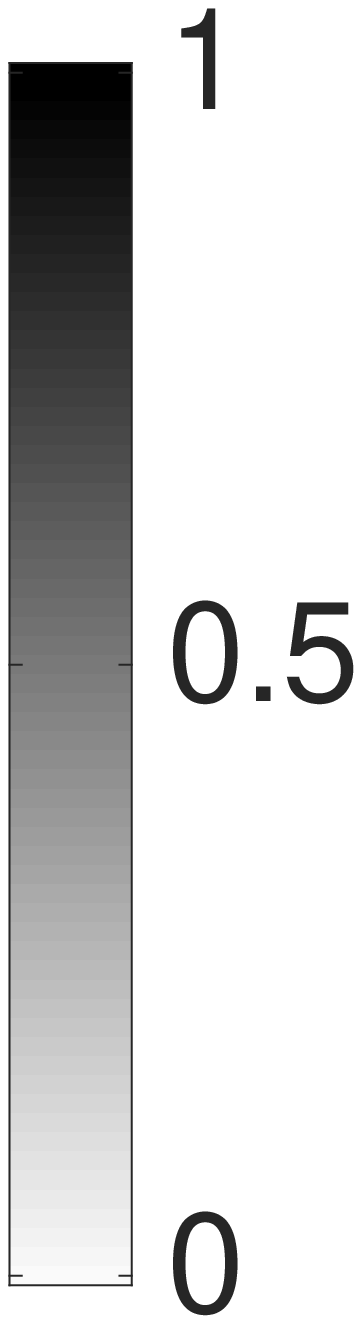}\\
\subfigure[1st coil]{\includegraphics[width=0.22\textwidth]{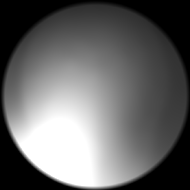}\label{subfig:gt2}}
\subfigure[2nd coil]{\includegraphics[width=0.22\textwidth]{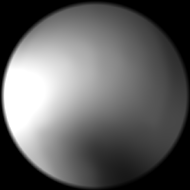}\label{subfig:gt3}}
\subfigure[3rd coil]{\includegraphics[width=0.22\textwidth]{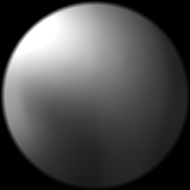}\label{subfig:gt4}}
\subfigure[4th coil]{\includegraphics[width=0.22\textwidth]{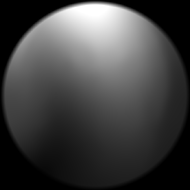}\label{subfig:gt5}}\hspace{0.1cm}\includegraphics[width=0.041\textwidth]{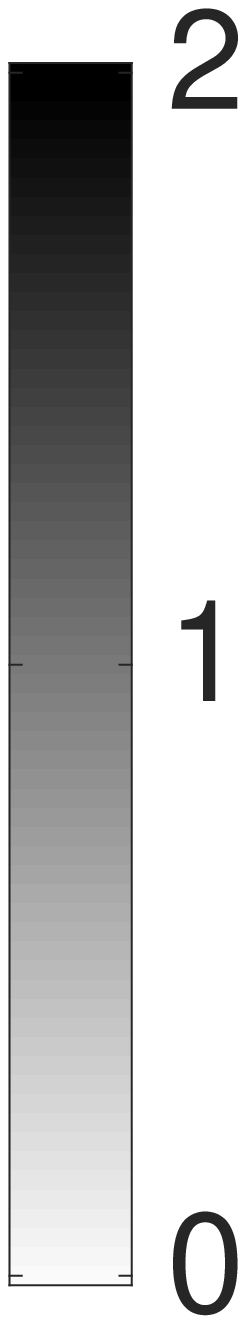}\\
\subfigure[5th coil]{\includegraphics[width=0.22\textwidth]{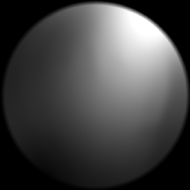}\label{subfig:gt6}}
\subfigure[6th coil]{\includegraphics[width=0.22\textwidth]{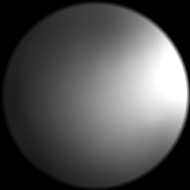}\label{subfig:gt7}}
\subfigure[7th coil]{\includegraphics[width=0.22\textwidth]{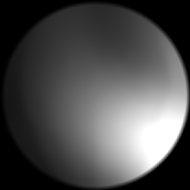}\label{subfig:gt8}}
\subfigure[8th coil]{\includegraphics[width=0.22\textwidth]{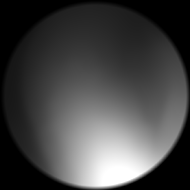}\label{subfig:gt5}\label{subfig:gt9}}\hspace{0.1cm}\includegraphics[width=0.041\textwidth]{Matlab/Reconstructions/LEGEND2.eps}
\end{center}
\caption{Figure \ref{subfig:gt1} shows the brain phantom as described in Section \ref{sec:expsetup}. Figure \ref{subfig:gt2} - \ref{subfig:gt9} show visualisations of the measured coil sensitivities of a water bottle. Figure \ref{subfig:sampscheme} shows the simulated, spiral-shaped sub-sampling scheme used to sub-sample the k-space data.}
\end{figure}

We want to demonstrate the numerical capabilities of the algorithm by applying it to the nonlinear problem of joint estimation of the spin-proton density and the coil sensitivities in parallel MRI. The discrete problem of joint reconstruction from sub-sampled k-space data on a rectangular grid reads
\begin{align*}
\left(\begin{array}{c} 
\hat{u}\\
\hat{c_1}\\
\vdots\\
\hat{c_2}
\end{array}\right)
\in \argmin_{\mathbf{v}=(u, c_1, \ldots, c_n)} \left\{ \frac{1}{2} \sum_{j = 1}^n \| S\mathcal{F}(G(\mathbf v))_j - f_j \|_2^2 + \alpha_0 R_0(u) + \sum_{j = 1}^n \alpha_j R_j(c_j) \right\} \, \text{,}
\end{align*} 
where $\mathcal{F}$ is the 2D discrete Fourier transform, $f_j$ are the k-space measurements for each of the $n$ coils, $S$ is the sub-sampling operator and $R_j$ denote appropriate regularisation functionals. The nonlinear operator $G$ maps the unknown spin-proton density $u$ and the different coil sensitivities $c_j$ as follows \cite{uecker2008image}:
\begin{align}
G(u, c_1, \ldots, c_n) = (u c_1, u c_2, \ldots, u c_n)^T \, \text{.}\label{eq:coilsensop}
\end{align}
In order to compensate for sub-sampling artefacts in sub-sampled MRI it is common practice to use total variation as a regulariser \cite{block2007undersampled,ramani2011parallel}. Coil sensitivities are assumed to be smooth, cf. Figure \ref{fig:expsetup}, motivating a reconstruction model similar to the one proposed in \cite{knoll2012parallel}. We therefore choose the discrete isotropic total variation, $R_0(u) = \|\nabla u \|_{2, 1}$, and the smooth 2-norm of the discretised gradient, i.e. $R_j(c_j) := \| \nabla c_j \|_{2, 2}$, for all $j > 0$, following the notation in \cite{benning2014phase}. We further introduce regularisation parameters $\lambda_j$ in front of the data fidelities and rescale all regularisation parameters such that $\alpha_0 + \sum_{j = 1}^n \lambda_j + \sum_{j = 1}^n \alpha_j = 1$. In order to realise this model via Algorithm \ref{alg:nlpadmm} we consider the following operator splitting strategy. We define $F(u_0, \ldots, u_n, v_0, \ldots, v_{2n})$ as
\begin{align*}
F(u_0, \ldots, u_n, v_1, \ldots, v_n) := \left( \begin{array}{c}
G(u_0, \ldots, u_n)\\ \begin{array}{cccc} \nabla u_0  & 0 & \cdots & 0\\ 0 & \nabla u_1 & \ddots & \vdots\\ \vdots & \ddots & \ddots & 0 \\ 0 & \cdots & 0 & \nabla u_n\end{array} 
\end{array}\right) - \left( \begin{array}{c}
v_0 \\ \vdots \\ v_n \\ \vdots \\ v_{2n}
\end{array}\right) \, \text{,}
\end{align*}  
set $H(u_0, \ldots, u_n) \equiv 0$, and $J(v_0, \ldots, v_{2n}) = \sum_{j = 0}^{2n} J_j(v_j)$ with $J_j(v_j) := \frac{\lambda_j}{2}\| S\mathcal{F} v_j - f_j \|_2^2$ for $j \in \{0, \ldots, n - 1\}$, $J_n(v_n) = \alpha_0 \| v_n \|_{2, 1}$ and $J_j(v_j) = \alpha_{j - n} \| v_j \|_{2, 2} $ for $j \in \{ n + 1, \ldots, 2n\}$. Note that with these choices of functions, all the resolvent operations can be carried out easily. In particular, we obtain
\begin{align*}
(I + \tau_1^k \partial H)^{-1}(w) &= w \, \text{,}\\
(I + \tau_2^k \partial J_j)^{-1}(w) &=  \mathcal{F}^{-1} \left( \frac{\mathcal{F}w_j + \tau_2^k \lambda_j S^T f_j}{1 + \tau_2^k \lambda_j \text{diag}(S^T S) } \right) \ \text{for} \ j \in \{0, \ldots, n - 1\} \, \text{,}\\
(I + \tau_2^k \partial J_n)^{-1}(w) &= \frac{w_n}{\| w_n \|_2}\max\left( \| w_n \|_2 - \alpha_0 \tau_2^k, 0 \right) \, \text{,}\\
(I + \tau_2^k \partial J_j)^{-1}(w) &= \frac{w_n}{\| w_n \|_{2, 2}}\max\left( \| w_n \|_{2, 2} - \alpha_j \tau_2^k, 0 \right) \ \text{for} \ j \in \{n + 1, \ldots, 2n\} \, \text{.}
\end{align*} 
Moreover, as $B_k = -I$ (and thus, $\| B^k \| = 1$) for all $k$, we can simply eliminate $\tau_2^k$ by replacing it with $1/\delta$, similar to
Section \ref{sec:nlpdhgm}. 

\begin{figure}[!ht]\label{fig:lownoise}
\begin{center}
\subfigure[Zero-filling]{\includegraphics[width=0.35\textwidth]{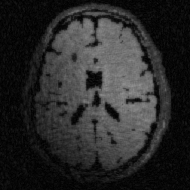}\label{subfig:gt1new}}
\hspace{0.1cm}
\subfigure[Reconstruction $u$]{\includegraphics[width=0.35\textwidth]{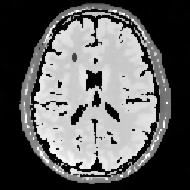}\label{subfig:recon1}}\hspace{0.25cm}\includegraphics[width=0.0975\textwidth]{Matlab/Reconstructions/LEGEND1.eps}\\
\subfigure[1st coil]{\includegraphics[width=0.22\textwidth]{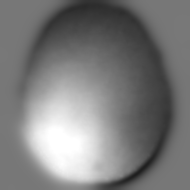}\label{subfig:recon2}}
\subfigure[2nd coil]{\includegraphics[width=0.22\textwidth]{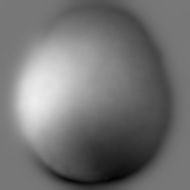}\label{subfig:recon3}}
\subfigure[3rd coil]{\includegraphics[width=0.22\textwidth]{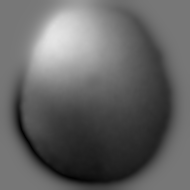}\label{subfig:recon4}}
\subfigure[4th coil]{\includegraphics[width=0.22\textwidth]{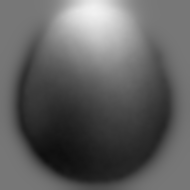}\label{subfig:recon5}}\hspace{0.1cm}\includegraphics[width=0.041\textwidth]{Matlab/Reconstructions/LEGEND2.eps}\\
\subfigure[5th coil]{\includegraphics[width=0.22\textwidth]{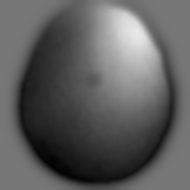}\label{subfig:recon2}}
\subfigure[6th coil]{\includegraphics[width=0.22\textwidth]{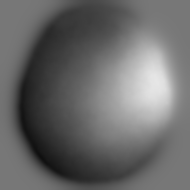}\label{subfig:recon3}}
\subfigure[7th coil]{\includegraphics[width=0.22\textwidth]{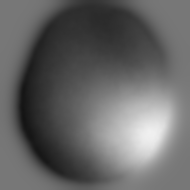}\label{subfig:recon4}}
\subfigure[8th coil]{\includegraphics[width=0.22\textwidth]{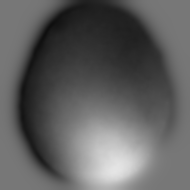}\label{subfig:recon5}}\hspace{0.1cm}\includegraphics[width=0.041\textwidth]{Matlab/Reconstructions/LEGEND2.eps}
\end{center}
\caption{Reconstructions for noise with low noise level $\sigma = 0.05$. Despite the sub-sampling, features of the brain phantom are very well preserved. In addition, the coil sensitivities seem to correspond well to the original ones, despite a slight loss of contrast. Note that coil sensitivities remain the initial value where the signal is zero.}
\end{figure}

\vspace{-0.25cm}
\subsection{Experimental setup}\label{sec:expsetup}
We now want to discuss the experimental setup. We want to reconstruct the synthetic brain phantom in Figure \ref{subfig:gt1} from sub-sampled k-space measurements. The numerical phantom is based on the design in \cite{aubert2006new} with a matrix size of $190 \times 190$. It consists of several different tissue types like cerebrospinal fluid (CSF), gray matter (GM), white matter (WM) and cortical bone. Each pixel is assigned a set of MR tissue properties: Relaxation times $\text{T}_1(x,y)$ and $\text{T}_2(x,y)$ and spin density $\rho(x,y)$. These parameters were also selected according to \cite{aubert2006new}. The MR signal $s(x,y)$ in each pixel was then calculated by using the signal equation of a fluid attenuation inversion recovery (FLAIR) sequence \cite{bernstein2004handbook}:
\begin{align*}
s(x,y) = \rho(x,y)(1-2 ~e^{-\text{TI}/\text{T}_1(x,y)})(1 -  e^{-\text{TR}/\text{T}_1(x,y)}) ~e^{-\text{TE}/\text{T}_2(x,y)}.
\end{align*}
The sequence parameters were selected: TR = 10000 ms, TE = 90 ms. TI was set to 1781 ms to achieve signal nulling of CSF ($\text{T}_1^\text{csf} \log(2)$ with $\text{T}_1^\text{csf} = 2569 \text{ms}$).

In order to generate artificial k-space measurements for each coil, we proceed as follows. First, we produce 8 images of the brain phantom multiplied by the measured coil sensitivity maps shown in Figure \ref{subfig:gt2} - \ref{subfig:gt9}. The coil sensitivity maps were generated from the measurements of a water bottle with an 8-channel head coil array. Then we produce artificial k-space data by applying the 2D discrete Fourier-transform to each of those individual images. Subsequently, we sub-sample only approx. 25\% of each of the k-space datasets via the spiral shown in Figure \ref{subfig:sampscheme}. Finally, we add Gau\ss ian noise with standard deviation $\sigma$ to the sub-sampled data.
 

\begin{figure}[!t]\label{fig:highnoise}
\begin{center}
\subfigure[Zero-filling]{\includegraphics[width=0.35\textwidth]{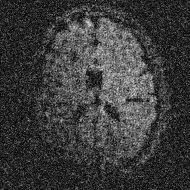}\label{subfig:gt1new}}
\hspace{0.1cm}
\subfigure[Reconstruction $u$]{\includegraphics[width=0.35\textwidth]{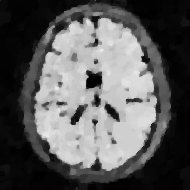}\label{subfig:recon1}}\hspace{0.25cm}\includegraphics[width=0.0975\textwidth]{Matlab/Reconstructions/LEGEND1.eps}\\
\subfigure[1st coil]{\includegraphics[width=0.22\textwidth]{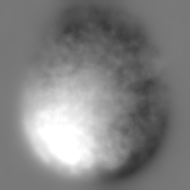}\label{subfig:recon2}}
\subfigure[2nd coil]{\includegraphics[width=0.22\textwidth]{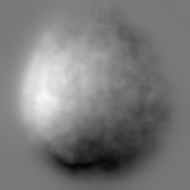}\label{subfig:recon3}}
\subfigure[3rd coil]{\includegraphics[width=0.22\textwidth]{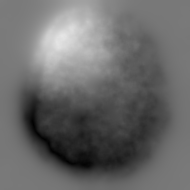}\label{subfig:recon4}}
\subfigure[4th coil]{\includegraphics[width=0.22\textwidth]{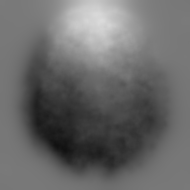}\label{subfig:recon5}}\hspace{0.1cm}\includegraphics[width=0.041\textwidth]{Matlab/Reconstructions/LEGEND2.eps}\\
\subfigure[5th coil]{\includegraphics[width=0.22\textwidth]{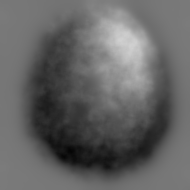}\label{subfig:recon2}}
\subfigure[6th coil]{\includegraphics[width=0.22\textwidth]{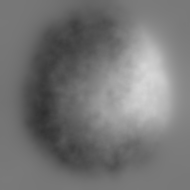}\label{subfig:recon3}}
\subfigure[7th coil]{\includegraphics[width=0.22\textwidth]{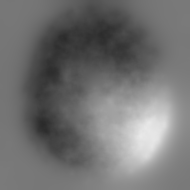}\label{subfig:recon4}}
\subfigure[8th coil]{\includegraphics[width=0.22\textwidth]{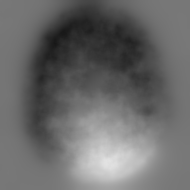}\label{subfig:recon5}}\hspace{0.1cm}\includegraphics[width=0.041\textwidth]{Matlab/Reconstructions/LEGEND2.eps}
\end{center}
\caption{Reconstructions for noise with high noise level $\sigma = 0.95$. Due to the large amount of noise, higher regularisation parameters are necessary. As a consequence, fine structures are smoothed out and in contrast to the case of little noise, compensation of sub-sampling artefacts is less successful.}
\end{figure}

\vspace{-0.25cm}
\subsection{Computations}
For the actual computations we use two noisy versions $f_j$ of the simulated k-space data; one with small noise ($\sigma = 0.05$) and one with a high amount of noise ($\sigma = 0.95$). As stopping criterion we simply choose a fixed number of iterations; for both the low noise level as well as the high noise level dataset we have fixed the number of iterations to 1500. The initial values used for the algorithm are $u_j^0 = \textbf{1}$ with $\textbf{1} \in \R^{l \times 1}$ being the constant one-vector, for all $j \in \{0, \ldots, n\}$. All other initial variables ($v^0$, $\mu^0$, $\overline{\mu}^0$) are set to zero. 

\subsubsection{Low noise level}
We have computed reconstructions from the noisy data with noise level $\sigma = 0.05$ via Algorithm \ref{alg:nlpadmm}, with regularisation parameters set to $\lambda_j = 0.0621$, $\alpha_0 = 0.062$ and $\alpha_j = 0.9317$ for $j \in \{1, \ldots, n\}$. We have further created a na\"{i}ve reconstruction by averaging the individual inverse Fourier-transformed images obtained from zero-filling the k-space data. The modulus images of the results are visualised in Figure \ref{fig:lownoise}. The results are visualised in Figure \ref{fig:highnoise}. The PSNR values for the averaged zero-filled reconstruction is 10.2185, whereas the PSNR of the reconstruction with the proposed method is 24.5572.

\subsubsection{High noise level}
We proceeded as in the previous section, but for noisy data with noise level $\sigma = 0.95$. The regularisation parameters were set to $\lambda_j = 0.0149$, $\alpha_0 = 0.0135$ and $\alpha_j = 0.9716$ for $j \in \{1, \ldots, n\}$. The modulus images of the results are visualised in Figure \ref{fig:highnoise}. The PSNR values for the averaged zero-filled reconstruction is 9.9621, whereas the PSNR of the reconstruction with the proposed method is 16.672.

\section{Conclusions \& outlook}
We have presented a novel algorithm that allows to compute minimisers of a sum of convex functionals with nonlinear operator constraint. We have shown the connection to the recently proposed NL-PDHGM algorithm which implies local convergence results in analogy to those derived in \cite{valkonen2014primal}. Subsequently we have demonstrated the computational capabilities of the algorithm by applying it to a nonlinear joint reconstruction problem in parallel MRI.

For future work, the convergence of the algorithm in the general setting has to be verified, and possible extensions to guarantee global convergence have to be studied. Generalisation of stopping criteria such as a linearised primal-dual gap will also be of interest as well. With respect to the presented parallel MRI application, exact conditions for the convergence (like the exact norm of the bounds) have to be verified. The impact of the algorithm- as well as the regularisation-parameters on the reconstruction has to be analysed, and  rigorous study with synthetic and real data would also be desirable. Moreover, future research focus will be on alternative regularisation functions, e.g. based on spherical harmonics motivated by \cite{sbrizzi2014robust}. Last but not least, other applications that can be modelled via \eqref{eq:nonlinprobform} should be considered in future research.

\subsubsection*{Acknowledgments} MB, CS and TV acknowledge EPSRC grant EP/M00483X/1. FK ackowledges National Institutes of Health grant NIH P41 EB017183.

\bibliography{nonlinearADMM}
\bibliographystyle{plain}

\end{document}